\documentclass[12pt]{article}
\usepackage{amssymb}
\usepackage{amsmath, amscd, amsfonts, array}
\begin{document}
\vspace{7cm}
\title{\bf{Some irreducibility results for truncated binomial expansions}}
\author{ Sudesh K. Khanduja\footnote{All correspondence
may be addressed to this author.},~ Ramneek Khassa and Shanta Laishram\\}

\date{}
\maketitle \baselineskip 18pt \noindent

\noindent\textbf{Abstract.} For positive integers $n>k$, let $P_{n,k}(x)=\displaystyle\sum_{j=0}^k \binom{n}{j}x^j $
be the polynomial obtained by truncating the binomial expansion of $(1+x)^n$ at the $k^{th}$ stage.
These polynomials arose in the investigation of Schubert calculus in Grassmannians.
In this paper, the authors prove the irreducibility of $P_{n,k}(x)$ over the field of rational numbers when
$2\leqslant 2k \leqslant n<(k+1)^3 $.\vspace{7mm}\vspace{2mm}
\vspace{5mm}\\
\noindent \textbf{1. Introduction}\vspace{3mm}\\
\indent For positive integers $k$ and $n$ with $k\leqslant n-1$, let $P_{n,k}(x)$
denote the polynomial $\displaystyle\sum_{j=0}^k \binom{n}{j}x^j$, where $\binom{n}{j}=\frac{n!}{j!~(n-j)!}$. In 2007, Filaseta, Kumchev
and Pasechnik considered the problem of irreducibility of $P_{n,k}(x)$ over the field $\mathbb{Q}$ of
rational numbers. This problem arose during the 2004 MSRI program on ``topological aspects
of real algebraic geometry" in the work of Inna Scherbak [6].  These
polynomials have also arisen in the context of work of Iossif V.
Ostrovskii [3]. In the case $k=2$, $P_{n,k}(x)$ has
negative discriminant and hence is irreducible over $\mathbb{Q}$. In fact it is already known that $P_{n,k}(x)$ is irreducible over $\mathbb{Q}$ for all $n\leqslant 100,~k+2\leqslant n$
(cf. [2, p.455]). In [2], Filaseta et al. pointed out that when $k=n-1$,
then $P_{n,k}(x)$ is irreducible over $\mathbb{Q}$ if and only if
$n$ is a prime number.
They also proved that for any fixed integer $k\geqslant 3$, there
exists an integer $n_0$ depending on $k$ such that $P_{n,k}(x)$ is
irreducible over $\mathbb{Q}$ for every $n\geqslant n_0.$ So there
are indications that $P_{n,k}(x)$ is irreducible for every
$n,k$ with $3\leqslant k\leqslant n-2$.\\
\indent In this paper, we prove the irreducibility of
$P_{n,k}(x)$ for all $n,k$ such that $2\leqslant 2k\leqslant
n<(k+1)^3$.  We consider the
irreducibility of the polynomial
$P_{n,k}(x-1)=\displaystyle\sum_{j=0}^k c_j x^j,$ where $c_j=\displaystyle\sum_{i=j}^k \binom{n}{i}\binom{i}{j} (-1)^{i-j}$. As in [2], on using the
identity\\
\centerline{$\displaystyle\sum_{j=0}^a (-1)^j\binom{b}{j}=(-1)^a\binom{b-1}{a}, \ a< b ~~\textrm{non-negative integers},$}\\
a simple calculation shows that
\begin{equation}
c_j=(-1)^{k-j} \binom{n}{j} \binom{n-j-1}{k-j}=\frac{(-1)^{k-j}~~ n(n-1)\cdots (n-k)}{j!(k-j)!}~\frac{1}{(n-j)}.
\end{equation}
In fact we shall prove the irreducibility of $P_{n,k}(x)$ using
Newton polygons with respect to primes exceeding $k$ dividing $\binom{n}{k}$
and some results of Erd\H{o}s, Selfridge, Saradha, Shorey and Laishram regarding such primes (cf. [7], [5]).
The same method gives the irreducibility of polynomial
\begin{equation} F_{n,k}(x)=\sum_{j=0}^k a_j c_j x^j, \end{equation}
where $a_0,a_1,\ldots,a_k$ are non-zero integers and each $a_i$ has all
of its prime factors  $\leqslant k$.\\  We prove\vspace{1mm}\\
\textbf{Theorem 1.1.} \emph{Let $k$ and $n$ be positive integers such that
$2k\leqslant n<(k+1)^3$. Then $P_{n,k}(x)$ is irreducible over $\mathbb{Q}$.}\vspace{3mm}\\
\indent Theorem 1.1 is derived from the following more general result.\vspace{3mm}\\
\textbf{Theorem 1.2.} \emph{Let $k$ and $n$ be positive integers such that $8\leqslant 2k \leqslant n<(k+1)^3$ and
$F_{n,k}(x)$ be as in (2). Then $F_{n,k}(x)$ is irreducible over $\mathbb{Q}$ except
possibly when $(n, k)$ belongs to the set} $ \{(8, 4), (10, 5), (12, 6), (16,8)\}.$\vspace{3mm}\\
\indent It may be pointed out that the polynomial\footnote{This
example was constructed by the referee.}
$F_{10,5}(x)$ given by \\
\centerline{$F_{10,5}(x)=2000.c_5x^5-375.c_4x^4-9.c_3x^3+10.c_2x^2-27.c_1x+25.c_0$}\\
\centerline{$=2000\cdot 252x^5+375\cdot 1050x^4-9\cdot 1800x^3-10\cdot 1575x^2-27\cdot 700x-25\cdot 126$}
\\has $7x^2+7x+1$ as a factor which shows that $F_{n,k}$ can be
reducible
over $\mathbb{Q}$.\vspace{3mm}\\ \indent In the course of the proof of Theorem 1.2, we prove the following result
 which is of independent interest as well.\newpage
\noindent \textbf{Theorem 1.3.} \emph{Let $k,n$ be integers such
that $n\geqslant k+2 \geqslant 4.$ Suppose there exists a prime
$p>k, p|(n-l)$ with $1\leqslant l\leqslant k-1$ and ord$_p(n-l)=e$
such that gcd$(e,l)\leqslant 2$ and gcd$(e,k-l)\leqslant 2$. If
$l_1< k/2$ is a positive integer such that $l \notin
\{l_1,2l_1,k-l_1,k-2l_1\}$, then $F_{n,k}(x)$ cannot have a factor
of degree $l_1$ over $\mathbb{Q}$.}\vspace{3mm}\\
\noindent\textbf{2. Notation and Preliminary Results}\vspace{3mm}\\
\indent For any non-zero integer $a$, let $v_p(a)$= ord$_p(a)$ denote the
$p$-adic valuation of $a$, i.e., the highest power of $p$ dividing
$a$ and denote $v_p(0)$ by $\infty$. Let
$g(x)=\displaystyle\sum_{j=0}^k a_j x^j$ be a polynomial over
$\mathbb{Q}$ with $a_0 a_k \neq 0$. To each term $a_i x^i$, we
associate a point $(n-i,v_p(a_i))$ ignoring however the point
$(n-i,\infty)$ if $a_i=0$ and form the set
$$S=\{(0,v_p(a_k)),\ldots,(n-j,v_p(a_j)),\ldots,(k,v_p(a_0))\}.$$
The Newton polygon of $g(x)$ with respect to $p$ is the polygonal path
formed by the lower edges along the convex hull of points of S.
Slopes of the edges are increasing when calculated from left to
right.\vspace{3mm}\\
\indent We begin with the following well known results (see [1]
for Theorem 2.A and
[4, 5.1.F] for Theorem 2.B).\vspace{4mm}\\
\textbf{Theorem 2.A.} \emph{Let $p$ be a prime and $g(x),h(x)$
belong to $\mathbb{Q}[x]$ with $g(0)h(0)\neq 0$ and $u\neq 0$ be
the leading coefficient of $g(x)h(x)$. Then the edges of the
Newton polygon of $g(x)h(x)$ with respect to $p$ can be formed by
constructing a polygonal path beginning at $(0,v_p(u))$ and using
the translates of the edges in the Newton polygon of $g(x)$ and
$h(x)$ with respect to $p$ taking exactly one translate for each
edge. The edges are translated in such a way as to form a
polygonal path with slopes of edges increasing.}
\vspace{4mm}\\
\textbf{Theorem 2.B.}
 \emph{Let $(x_0,y_0),(x_1,y_1),\ldots,
(x_r,y_r)$ denote the successive vertices of the Newton polygon of
a polynomial $g(x)$ with respect to a prime $p$. Let $\tilde{v}_p$
denote the unique extension of ${v}_p$ to the algebraic closure of
$\mathbb{Q}_p$, the field of $p$-adic numbers. Then $g(x)$ factors over
$\mathbb{Q}_p$ as $g_1(x)g_2(x) \cdots g_r(x)$ where the
degree of $g_i(x)$ is $x_i-x_{i-1},~i=1,2,\ldots,r$ and all the roots
of $g_i(x)$ in the algebraic closure of $\mathbb{Q}_p$ have
$\tilde{v}_p$ valuation $\frac{y_i-y_{i-1}}{x_i-x_{i-1}}$. In
particular all the roots of an irreducible factor of $g(x)$ over
$\mathbb{Q}_p$ will have the same
$\tilde{v}_p$ valuation.}
\vspace{4mm}\\
\indent  For an integer $\nu>1$, let $P(\nu)$ denote the greatest
prime divisor of $\nu$ and let  $\pi(\nu)$ denote the number of
primes not exceeding $\nu$. As in [5], $\delta(k)$ will denote the
integer\newpage \noindent defined for $k\geqslant 3 $ by \vspace{2mm}\\
$\delta(k)= \left\{
              \begin{array}{ll}
                2, & \hbox{if~ $3\leqslant k\leqslant 6$;} \\
                1, & \hbox{if~$7\leqslant k\leqslant 16$;} \\
                0, & \hbox{otherwise.}
              \end{array}
            \right.$\vspace{2mm}\\
 For numbers $n, k$ and $h$, $[n, k, h]$ will stand for the set of all pairs $(n, k), (n+1, k), \ldots$ , $(n+h-1, k)$. In particular
$[n,k,1]=\{(n, k)\}$.\vspace{2mm}\\
\indent We shall denote  by $S$ the union of
the sets\vspace{4mm}\\
\noindent$[6,3,1],[8,3,3],[18,3,1],[9,4,1],[10,5,4],[16,5,1],[18,5,3],[27,5,2],[12,6,2],[20,6,1],$\\ $[14,7,3],[18,7,1],[20,7,2],[30,7,1],[16,8,1],[21,8,1],[26,13,3],[30,13,1],[32,13,2],$ \\ $[36,13,1],[28,14,1],[33,14,1],[36,17,1]$ \vspace{2mm}\\ and by $T$ the union of the sets \vspace{2mm}\\
$[38,19,3],[42,19,1],[40,20,1],[94,47,3],[100,47,1],[96,48,1],[144,71,2],[145,72,1],$\\$[146,73,3], [156,73,1],[148,74,1],[162,79,1],[166,83,1],[172,83,1],[190,83,1],$ \\ $[192,83,1],[178,89,1],[190,89,1],[192,89,1],[210,103,2],[212,103,2] [216,103,2],$\\$[213,104,1],[217,104,1],[214,107,12],[216,108,10],[218,109,9], [220,110,7]$\\$[222,111,5],[224,112,3],[226,113,7],[250,113,1],[252,113,2],[228,114,5],
[253,114,1],$\\$[230,115,3],[232,116,1],[346,173,1],[378,181,1],[380,181,2],[381,182,1],
[392,193,2],$\\$[393,194,1],[396,197,1],[398,199,3],[400,200,1],[552,271,5],[553,272,1],
[555,272,2],$\\$[556,273,1],[554,277,3],[558,277,5],[556,278,1],[559,278,4],[560,279,3],
[561,280,1],$\\$[562,281,7],[564,282,5],[566,283,5],[576,283,1],[568,284,3],[570,285,1],
[586,293,1]$.\vspace{4mm}\\
\indent With the above notations, we shall use the following theorem due to Laishram and Shorey $\textrm{[5, Theorem 3]}$.\vspace{3mm}\\
\noindent\textbf{Theorem 2.C.}
\emph{Let $n\geqslant 2k\geqslant 6$ and $f_1<f_2<\cdots<f_\mu$ be integers in $[0, k)$. Assume that the
greatest prime factor of $(n-f_1)\ldots(n-f_\mu)\leqslant k$. Then
$\mu\leqslant k-\left[\frac{3}{4} \pi(k) \right] +1-\delta(k)$ unless $(n,k) \in S \cup T$.}\vspace{2mm}\\
\indent The following corollary is an immediate consequence of Theorem 2.C.\vspace{3mm}\\
\noindent\textbf{Corollary 2.D.} \emph{Let $n$ and $k$ be positive integers with $n\geqslant 2k \geqslant 38$.
Then there are at least five distinct terms of the product $n(n-1) \cdots (n-k+1)$
each divisible by a prime exceeding $k$ except when $(n, k) \in T$.}\vspace{3mm}\\
\indent For the proof of Theorem 1.3, we need the following propositions.\vspace{3mm}\\
\textbf{Proposition 2.1.}\emph{ Let $k\geqslant 6$ and $n > k^2$. Then there exist two distinct terms
$n+r$ and $n+s$ of the product $n(n + 1)\cdots (n+k-1)$
which are divisible by primes $> k$ exactly to an odd power.}\\
\emph{Proof.} Suppose the proposition is false for some $n$ and $k$ with $k\geqslant 6$ and $n>k^2$. Let
$\Delta(n, k)=n(n + 1)\cdots (n+k-1)$. Thus either ord$_p(\Delta(n, k))$ is even for all primes $p>k$
or there is exactly one term $n+i$ and a prime $p>k$ such that ord$_p(\Delta(n, k))$ is odd.
The first possibility is excluded since for any positive integer $b$ with
$P(b)\leqslant k $, the equation $$ n(n + 1)\cdots (n+k-1)=by^2 $$
has no solution in positive integers $n,k,y$ when $n> k^2\geqslant 4^2 $ by [7, Theorem A].
We now consider the case when there is exactly a term $n+i$ and a prime $p>k$ such that
ord$_p(\Delta(n, k))$ is odd. Suppose first that $0< i< k-1$. Removing the term $n+i$ from $\Delta (n,k)$,
we see that $n(n+1)\cdots (n+i-1)(n+i+1)\cdots (n+k-1)=b_1 y_1^2$ where $P(b_1)\leqslant k$ which is
impossible by virtue of [7, Theorem 2\footnote{It states that for  $n>k^2\geqslant 5^2$  the equation
$n(n+1)\cdots (n+i-1)(n+i+1)\cdots (n+k-1)=by^2$
has no solution in positive integers $n, k, b, y$ with $P(b)\leqslant k$ and $0<i<k-1$.}].

It remains to consider the case when $i=0$ or $k-1$. Let $\Delta'$ denote the product
$(n + 1)\cdots (n+k-1)$ or $n(n + 1)\cdots (n+k-2)$ according as $i=0$ or $k-1$.
Then $\Delta'$ is a product of $k-1$ consecutive integers such that
\begin{equation} \Delta'=b_2 y_2^2 \end{equation}
with $P(b_2)\leqslant k$. This is impossible when $P(b_2)\leqslant
k-1$ by [7, Theorem A]. It only remains to deal with the situation
when $P(b_2)=k$. Then $k$ will be a prime dividing only one term
of the product $\Delta'$, say $k$ divides $n+j,j\neq i$. We remove
the term $n+j$ of the product $\Delta'$ and it is clear from (3)
that
\begin{equation}
\frac{\Delta'}{n+j}=b_3 y_3^2~,~~P(b_3)\leqslant k-2.
\end{equation}
It is immediate from (4) and [7, Theorem 2] that $n+j$ is the first or last term of the product $\Delta'$ as
$k-1 \geqslant 5.$ Thus we see that $ \frac{\Delta'}{n+j}$ is the product of $k-2$ consecutive integers.
This is impossible by [7, Theorem A].\vspace{3mm}\\
\textbf{Proposition 2.2.}\emph{ Let $n,k$ be positive integers with $n
\geqslant k+2 \geqslant 4$ and $F_{n,k}(x)$ be given by (2).
Suppose there exists a prime $p>k$ such that $p^e||(n-l)$ for some
$l,~1\leqslant l \leqslant k-1$. Let $d=$gcd$(e, l)$ and $d'=$gcd$(e, k-l)$.
Then the following hold.\\
(i) The edges of the Newton polygon of $F_{n,k}(x)$ with respect
to $p$ have slopes $\frac{-e}{k-l},~\frac{e}{l}$.\\
(ii) $F_{n,k}(x)$ has at least two distinct irreducible factors
over $\mathbb{Q}_p$; one of them has degree a multiple of
$\frac{l}{d}$ and other has degree a multiple of
$\frac{k-l}{d'}$.\\
(iii) If $d=d'=1$, then
$F_{n,k}(x)$ factors over $\mathbb{Q}_p$ as a product of two distinct irreducible
polynomials of degrees $l$ and $k-l$.}\vspace{1mm}\\
\emph{Proof.} We consider the Newton polygon of $F_{n,k}(x)$
with respect to the prime $p$. In view of (1), the vertices of the Newton polygon
are $(0,e),(k-l,0),(k,e) $. Thus the Newton polygon has two edges,
one from $(0,e)$ to $(k-l,0)$ and other from $(k-l,0)$ to
$(k,e)$ with respective slopes $\frac{-e}{k-l}$ and
$\frac{e}{l}$ proving $(i)$.

Note that equations of the two edges are given by:
$$y-e=\frac{-e}{k-l}~x  \ \ \textrm{and} \ \ y=\frac{e}{l}~(x-k+l).$$
On the first edge, the $x$-coordinates of the lattice points
occur at multiples of $\frac{k-l}{d'}$, i.e., when
$x=\frac{k-l}{d'}.M$ where $0\leqslant M\leqslant d'$;
on the second edge the $x$-coordinates of lattice points are given
by $k-l+\frac{l}{d}.N$ where $0\leqslant N \leqslant d$.
By Theorem 2.B, all the roots of an irreducible factor of $F_{n,k}(x)$
over $\mathbb{Q}_p$ have the same slope. Since the slopes of the two edges as shown in $(i)$ are
different, we see that any irreducible factor of $F_{n,k}(x)$ over $\mathbb{Q}_p$ must lie
on the first edge or on the second edge. Hence assertion $(ii)$ now follows from Theorem 2.A.
Assertion $(iii)$ is an immediate consequence of $(ii)$. The last assertion quickly yields the following result.\vspace{4mm}\\
\textbf{Corollary 2.3. } \emph{If for a pair $(n,k),~n\geqslant k+2$, there exist terms $n-l',~n-l'',~1\leqslant l'<l''<k$, divisible respectively by primes $p',p''$ exceeding $k$ exactly to the first power such that $l'+l''\neq k$, then $F_{n,k}(x)$ is irreducible over $\mathbb{Q}$.}\vspace{3mm}\\
\indent The following proposition is already known (cf. [2, Lemma 1]). For the sake of reader's convenience, it is proved here.\vspace{3mm}\\
\noindent \textbf{Proposition 2.4.} \emph{Let $n,k$ and $F_{n,k}(x)$ be as in Proposition 2.2. Let $p$ be a prime $>k$
and $e>0$ be such that $p^e||n$. Then every irreducible factor of
$F_{n,k}(x)$ over $\mathbb{Q}_p$ has degree a multiple of $\frac{k}{D}$, where
$D=$gcd$(e, k)$.}\vspace{1mm}\\
\emph{Proof.} The vertices of the Newton polygon of $F_{n,k}(x)$ with respect to $p$
are $(0,e),(k,0).$ Thus the Newton polygon has only one edge whose equation
is given by $y-e=\frac{-e}{k}~x.$ The $x$-coordinates of the lattice points
on this edge occur at multiples of $k/D$. So arguing as in Proposition 2.2,
any irreducible factor of $F_{n,k}(x)$ must have degree a multiple of $k/D.$\newpage
\noindent \textbf{3. Proof of Theorem 1.3}\vspace{3mm}\\
\indent As pointed out in the proof of Proposition 2.2 (with $d,~d'$ atmost 2), if $(x,y)$ is a lattice point on the
 Newton polygon of $F_{n,k}(x)$ with respect to $p$, then $x\in X=\{0,\frac{k-l}{2},k-l,k-\frac{l}{2},k\}.$
  By Theorems 2.A, 2.B, each irreducible factor of $F_{n,k}(x)$ over $\mathbb{Q}$ must have degree equal to a sum of numbers (may be one of the numbers) from \begin{align*} \frac{l}{2},\frac{l}{2},\frac{k-l}{2},\frac{k-l}{2};\end{align*}
these correspond to possible differences $x_i-x_{i-1}$ in Theorem 2.B, with the actual differences possibly
formed from sums of these possible differences. Thus an irreducible factor of $F_{n,k}(x)$ over $\mathbb{Q}$ must
 have degrees in the set $$\left\{\frac{l}{2},l,\frac{k}{2},\frac{k-l}{2},k-l,\frac{2k-l}{2},\frac{k+l}{2},k\right\}.$$
Given that $l<k$, the elements of this set that can be less than $k/2$ are $l/2,l,(k-l)/2$ and $k-l$.
The conditions in Theorem 1.3 imply that $l_1$ is not among $l/2,l,(k-l)/2$ and $k-l$, so the theorem follows.\vspace{3mm}\\
\textbf{4. Proof of Theorem 1.2}\vspace{3mm}\\
\indent With $S$ and $T$ as in Theorem 2.C, we first prove \vspace{2mm}\\
\textbf{Lemma 4.1.} \emph{For $(n,k)\in S\cup T,~k\geqslant 4,$
$F_{n, k}(x)$ is irreducible over $\mathbb{Q}$ except possibly
when $(n,k)$ belongs to the subset $S'$ of $S$ given by
$S'=\{(10,5),(12,6),(16,8)\}$.}\\
Proof. Let $S''$ denote the subset of $S$ given by
$S''=\{(9,4),(12,5),(16,5),(18,5),(27,5)\}$. Observe that if $n$
is divisible by a prime $p>k$ with ord$_p(n)=1$, then $x^k F_{n,
k}(1/x)$ is an Eisenstein polynomial with respect to $p$ and so
$F_{n,k}(x)$ is irreducible over $\mathbb{Q}$. Further if two
distinct terms $n-l_1,n-l_2$ of the product $n(n-1)\cdots (n-k+1)$
are divisible by primes $p_1$ and $p_2$ exceeding $k$ such that
ord$_{p_i}(n-l_i)=1$ and $l_1+l_2\neq k$, then in view of the
above observation and Corollary 2.3, $F_{n, k}(x)$ is irreducible
over $\mathbb{Q}$. For each $(n,k)$ belonging to $ T\cup
(S\setminus S'\cup S'')$ with  $n$ not divisible by any prime $>k$ up to
the first power, Table 1 at the end of this section indicates two
primes $p_1$ and $p_2$ satisfying the above property. It can be
easily seen that for $(n,k)\in S''$, $F_{9,4}(x)$ is an Eisenstein
polynomial with respect to the prime $5$, $F_{12,5}(x)$ is
Eisenstein with respect to $7$, $F_{16,5}(x),F_{27,5}(x)$ are
Eisenstein with respect to $11$ and  $F_{18,5}(x)$ is Eisenstein with
respect to $13$. Hence the
lemma is proved.\newpage
\noindent \textbf{Lemma 4.2.} \emph{For $8\leqslant n <5^3,$ the polynomial
$F_{n,4}(x)$ is irreducible over $\mathbb{Q}$ except when $n$
belongs to the set $U=\{8,50,98,100\}$}.\\
\emph{Proof}. As pointed out in the proof of Lemma 4.1, we need to
verify the irreducibility of $F_{n,4}(x)$ when $n$ is not
divisible by any prime more than $4$ exactly with the first power.
For such $n$ not exceeding $124$ and $n$ not belonging to the set $\{8,9,18,27,50,98,100\}$, Table 2 at the end of this
section indicates two terms $n-l',n-l'',~1\leqslant
l',l''\leqslant 3,~l'+l''\neq 4$ such that $n-l',n-l''$ are
divisible by primes $p',p''$ (respectively) up to the first power
only. So the lemma is proved in view of Corollary 2.3 and the fact
that $F_{9,4}(x),~F_{18,4}(x)$ and $F_{27,4}(x)$ are Eisenstein polynomials with respect
to the primes $5,7$ and $23$
respectively.\vspace{3mm}\\
 \emph{Proof of Theorem 1.2.} We divide the proof into
two cases. \\\textbf{Case I}.  ~$8\leqslant 2k\leqslant n<(k+1)^2$.
Note that the theorem is already proved in the present case for $k=4$
by virtue of Lemma 4.2, so it may be assumed that $k\geqslant 5$
here. Applying Theorem 2.C, we see that there exist at least three
terms $n-l_i,~i\in \{1,2,3\}$ which are divisible by primes
exceeding $k$ exactly up to the first power unless $(n,k)\in S\cup
T.$ Using Proposition 2.2 $(iii)$, $F_{n,k}(x)$ factors over
$\mathbb{Q}_{p_i}$ as a product of two non-associate irreducible
polynomials of degree $l_i$ and $k-l_i$ for $1 \leqslant i
\leqslant 3$. If $F_{n,k}(x)$ were reducible over $\mathbb{Q}$,
then $F_{n,k}(x)$ will have a factorization of the type
$F_{n,k}(x)=a_kc_kG_i(x) H_i(x)$ where $G_i(x),H_i(x)$ are monic
irreducible polynomials belonging to $\mathbb{Q}[x]$ with degrees
$k-l_i$, $l_i$ respectively. This is impossible as $l_1,l_2$ and
$l_3$ are distinct. So the theorem is proved in the present case
when $(n,k)$ does not belong to $S\cup T.$ When $(n,k)\in
(S\setminus S') \cup T$ with $k\geqslant 4$, the irreducibility of $F_{n,
k}(x)$ follows from
Lemma 4.1.\vspace{2mm}\\
\textbf{Case II}. $k\geqslant 4,~(k+1)^2\leqslant n <(k+1)^3.$ In
this case, we first show that $F_{n,k}(x)$ cannot factor over
$\mathbb{Q}$ as a product of two irreducible polynomials of degree
$\frac{k}{2}$ each. For this it is enough to show that there
exists $l'\neq k/2,~0\leqslant l' \leqslant k-1$ such that $n-l'$
is divisible by a prime $p'>k$ exactly with the first power. If
$l'=0$, then as pointed out in the opening lines of the proof of
Lemma 4.1, $F_{n, k}(x)$ is irreducible over $\mathbb{Q}$. If
$l'\geqslant 1$ then by Proposition 2.2 $(iii)$, $F_{n,k}(x)$ has
two irreducible factors of degree $l'$ and $k-l'$ over
$\mathbb{Q}_{p'}$. This leads to a contradiction as $l'\neq k/2$
thereby proving the irreducibility of $F_{n, k}(x)$ over
$\mathbb{Q}$. The existence of a term $n-l'\neq
n-\frac{k}{2},~0\leqslant l' \leqslant k-1,$ which is divisible by some
prime $p'>k$ with ord$_{p'}(n-l')=1$ is guaranteed for $k\geqslant
6$ by Proposition
2.1 as $(k+1)^2\leqslant n <(k+1)^3$ in the present situation. This proves the assertion stated in the opening lines of Case II.\\
\indent It only remains to be shown that $F_{n,k}(x)$ cannot
have a factor of degree less than $\frac{k}{2}$ over $\mathbb{Q}$. Suppose
to the contrary that it has a factor of degree
 $l_1<\frac{k}{2}$ over $\mathbb{Q}$. We make some claims.

\noindent
{\bf Claim 1:} $P(n)\leqslant k$. \\
Suppose not. Let $p$ be a prime $>k$ dividing $n$ with exact power $e\geqslant 1.$ Then
$e\leqslant 2$ since $n<(k+1)^3$. So by Proposition 2.4, every
irreducible factor of $F_{n,k}(x)$ over $\mathbb{Q}_{p}$ has
degree a multiple of $k$ or $\frac{k}{2}$ according as $e =1$ or
$2$ respectively. This is not possible in view of our supposition.

\noindent {\bf Claim 2:} There are at most four distinct terms in
the product
$n(n-1)\cdots (n-k+1)$ each of which is divisible by some prime $>k$. \\
Assume the contrary. Then there is a term $n-l$ with $0\leqslant
l<k$ and a prime $p>k$ with $p$ dividing $(n-l)$ such that
$l\notin \{l_1,~2l_1,~k-l_1,~k-2l_1 \}$ where $l_1$ is as in the
paragraph preceeding Claim I. Note that $l>0$ in view of Claim 1.
Further $e=$ord$_p(n-l)\leqslant 2$ implying that $F_{n,k}(x)$
cannot have a factor of degree $l_1$ over $\mathbb{Q}$ by Theorem
1.3, which contradicts our assumption.

 \noindent {\bf Claim 3:}
There are at most two distinct terms in the product
$n(n-1)\cdots (n-k+1)$ which are divisible by a prime $>\sqrt{n}$. \\
Suppose not. Let $1\leqslant l'_1<l'_2<l'_3~$ be such that there
exist primes $p_i>\sqrt{n}$ dividing $n-l'_i$. Note that ord$_{p_i}(n-l'_i)=1$ for $i\in \{1, 2, 3\}$. Since $(k+1)^2\leqslant n$,
in view of Proposition 2.2 $(iii)$, it follows that
$F_{n,k}(x)$ factors over $\mathbb{Q}_{p_i}$ as a product of two
non-associate irreducible polynomials of degree $l'_i$ and
$k-l'_i$, $1 \leqslant i \leqslant 3$. Arguing as in Case I, we get a contradiction because
$l'_1,l'_2$ and $l'_3$ are distinct.

From Claim 2, Corollary 2.D and Lemma 4.1, it follows that
$k\leqslant 18$. Note that for $k=4,$ in view of Lemma 4.2, we have
only to consider $n=50,98,100$ as $5^2\leqslant n< 125$. For each of
these values of $n$, $F_{n,k}(x)$ must be irreducible over
$\mathbb{Q}$ by virtue of Claim 1, as $P(n)$ is more than $4.$ For
$k\geqslant 5$, by virtue of Claim 1, we may first
 restrict to those $n$ for which $P(n)\leqslant k$. Further by Claims 2 and 3, those $n$ can be excluded for which
$n(n-1)\cdots (n-k+1)$ has either five terms divisible by a prime
$>k$ or three terms divisible by a prime $>\sqrt{n}$. We use
\emph{Sage} mathematics software for the above computations. Then
we are left with the following pairs $(n, k)$ given by
$$(50,5),(64,5),(100,5),(128,5),(200,5),(50,6).$$
All these pairs satisfy the hypothesis of
Corollary 2.3 as is clear from Table 3. This completes the proof of the
theorem.\vspace{3mm}\\
\noindent \textbf{Table 1.}
$$\begin{array}{|l|l|l|} \hline  (n,k)\in [n,k,h] \rightarrow Primes  & (n,k)\in [n,k,h]  \rightarrow   Primes  &  (n,k)\in [n,k,h]  \rightarrow  Primes  \\
\hline   [20,5,1]~~~~~~~~~~~~~~~17,19 &[162,79,1]~~~~~~~~~~~~~131,139   &[346,173,1]~~~~~~~~~~~~293,307  \\
\hline   [20,6,1]~~~~~~~~~~~~~~~17,19&[166,83,1]~~~~~~~~~~~~~131,139  &[378,181,1]~~~~~~~~~~~~293,307  \\
\hline   [14,7,3]~~~~~~~~~~~~~~~11,13&[172,83,1]~~~~~~~~~~~~~137,139 &[380,181,2]~~~~~~~~~~~~293,307  \\
\hline   [18,7,1]~~~~~~~~~~~~~~~13,17 &[190,83,1]~~~~~~~~~~~~~131,139  &[381,182,1]~~~~~~~~~~~~293,307  \\
\hline   [20,7,1]~~~~~~~~~~~~~~~17,19 &[192,83,1]~~~~~~~~~~~~~131,139   &[392,193,2]~~~~~~~~~~~~293,307  \\
\hline   [21,7,1]~~~~~~~~~~~~~~~17,19&[178,89,1]~~~~~~~~~~~~~131,139  &[393,194,1]~~~~~~~~~~~~293,307  \\
\hline   [30,7,1]~~~~~~~~~~~~~~~13,29&[190,89,1]~~~~~~~~~~~~~131,139  &[396,197,1]~~~~~~~~~~~~293,307  \\
\hline   [21,8,1]~~~~~~~~~~~~~~~17,19&[192,89,1]~~~~~~~~~~~~~139,149  &[398,199,3]~~~~~~~~~~~~293,307  \\
\hline   [26,13,3]~~~~~~~~~~~~~~19,23&[210,103,1]~~~~~~~~~~~~139,149   &[400,200,1]~~~~~~~~~~~~283,307  \\
\hline   [30,13,1]~~~~~~~~~~~~~~19,23&[212,103,2]~~~~~~~~~~~~139,149 &[552,271,5]~~~~~~~~~~~~421,431  \\
\hline   [32,13,2]~~~~~~~~~~~~~~29,31&[216,103,2]~~~~~~~~~~~~139,149  &[553,272,1]~~~~~~~~~~~~421,431  \\
\hline   [36,13,1]~~~~~~~~~~~~~~29,31&[213,104,1]~~~~~~~~~~~~139,149  &[555,272,2]~~~~~~~~~~~~421,431  \\
\hline   [28,14,1]~~~~~~~~~~~~~~17,19&[217,104,1]~~~~~~~~~~~~139,149   &[556,273,1]~~~~~~~~~~~~421,431  \\
\hline   [33,14,1]~~~~~~~~~~~~~~29,31&[214,107,12]~~~~~~~~~~~139,149  &[554,277,3]~~~~~~~~~~~~421,431  \\
\hline   [36,17,1]~~~~~~~~~~~~~~29,31&[216,108,10]~~~~~~~~~~~139,149 &[558,277,5]~~~~~~~~~~~~421,431  \\
\hline   [38,19,3]~~~~~~~~~~~~~~23,29&[218,109,9]~~~~~~~~~~~~139,149   &[556,278,1]~~~~~~~~~~~~421,431  \\
\hline   [42,19,1]~~~~~~~~~~~~~~37,41&[220,110,7]~~~~~~~~~~~~139,149 &[559,278,4]~~~~~~~~~~~~421,431  \\
\hline   [40,20,1]~~~~~~~~~~~~~~31,37&[222,111,5]~~~~~~~~~~~~139,149  &[560,279,3]~~~~~~~~~~~~421,431  \\
\hline   [94,47,3]~~~~~~~~~~~~~~89,83&[224,112,3]~~~~~~~~~~~~139,149   &[561,280,1]~~~~~~~~~~~~421,431  \\
\hline   [100,47,1]~~~~~~~~~~~~~83,89&[226,113,7]~~~~~~~~~~~~139,149  &[562,281,7]~~~~~~~~~~~~409,431  \\
\hline   [96,48,1]~~~~~~~~~~~~~~79,83&[250,113,1]~~~~~~~~~~~~139,149  &[564,282,5]~~~~~~~~~~~~409,431  \\
\hline   [144,71,2]~~~~~~~~~~~~~101,103&[252,113,2]~~~~~~~~~~~~139,149  &[566,283,5]~~~~~~~~~~~~421,431  \\
\hline   [145,72,1]~~~~~~~~~~~~~101,103&[228,114,5]~~~~~~~~~~~~139,149 &[576,283,1]~~~~~~~~~~~~421,431  \\
\hline   [146,73,3]~~~~~~~~~~~~~101,103&[253,114,1]~~~~~~~~~~~~139,149 &[568,284,3]~~~~~~~~~~~~419,431  \\
\hline   [156,73,1]~~~~~~~~~~~~~109,113&[230,115,3]~~~~~~~~~~~~139,149  &[570,285,1]~~~~~~~~~~~~421,431  \\
\hline [148,74,1]~~~~~~~~~~~~~107,113&[232,116,1]~~~~~~~~~~~~139,149  &[586,293,1]~~~~~~~~~~~~421,431  \\
 \hline
\end{array}$$
\textbf{Table 2.}\\
$$\begin{array}{|l|l|l|} \hline  n~ \rightarrow~ n-l',n-l'',~p',p''~   & n~ \rightarrow~~ n-l',n-l'',~p',p''~  &  n~ \rightarrow~~~ n-l',n-l'',~p',p''  \\
\hline 12~~~~~~~~~~~10,11,~~~~5,11  &48~~~~~~~~~~~~46,47,~~~~23,47  &81~~~~~~~~~~~~~79,80,~~~~~79,5   \\
\hline 16~~~~~~~~~~~14,15,~~~~7,5  &49~~~~~~~~~~~~46,47,~~~~23,47 &96~~~~~~~~~~~~~94,95,~~~~~47,19 \\
\hline 24~~~~~~~~~~~22,23,~~~~11,23  &54~~~~~~~~~~~~52,53,~~~~13,53  &108~~~~~~~~~~~106,107,~~~53,107  \\
\hline 25~~~~~~~~~~~22,23,~~~~11,23  &64~~~~~~~~~~~~62,63,~~~~31,7  &121~~~~~~~~~~~119,120,~~~17,5  \\
\hline 32~~~~~~~~~~~30,31,~~~~5,31  &72~~~~~~~~~~~~70,71,~~~~5,71 &  \\
\hline 36~~~~~~~~~~~34,35,~~~~17,5  &75~~~~~~~~~~~~73,74,~~~~73,37  &\\
\hline
\end{array}$$\newpage
\noindent\textbf{Table 3.}\\
$$\begin{array}{|l|l|l|} \hline  (n,k)~~~ \rightarrow~~~ n-l',n-l''   & (n,k)~~~ \rightarrow~~~ n-l',n-l''  &  (n,k)~~~~ \rightarrow~~~~ n-l',n-l''  \\
\hline (50,5)~~~~~~~~~~~~~~46,47  &(100,5)~~~~~~~~~~~~~97,99   &(200,5)~~~~~~~~~~~~~~197,199  \\
\hline (64,5)~~~~~~~~~~~~~~61,63  &(128,5)~~~~~~~~~~~~~126,127  &(50,6)~~~~~~~~~~~~~~~~46,47 \\
\hline
\end{array}$$
\vspace{2mm}\\
\noindent\textbf{5. Proof of Theorem 1.1}\vspace{3mm}\\
\indent In view of Theorem 1.2., we need to prove the irreducibility of $P_{n,k}(x)$ only when $1\leqslant k\leqslant 3$ with $2k\leqslant n<(k+1)^3$ or $(n, k)$ belongs to $\{(8, 4), (10, 5), (12, 6), (16,8)\}.$ Using Maple, we have verified the irreducibility of $P_{n,k}(x)$ for these values of $(n, k)$.\vspace{3mm}\\
\noindent\textbf{Acknowledgements}\vspace{2mm}\\
\indent The financial support by National Board for Higher
Mathematics, Mumbai and by CSIR (grant no. 09/135(0525)/2007-EMR-I) is gratefully
acknowledged. The authors are thankful to the referee for several helpful suggestions.\vspace{3mm}\\
\noindent \textbf{References}.\vspace{3mm}\\
{[1] M. G. Dumas, Sur quelques cas d'irr\'{e}ducibilit\'{e} des polynomes \`{a} coefficients rationnels, \textit{J. Math. Pures Appl.} 2 (1906) 191-258.\\
{[2] M. Filaseta, A. Kumchev and D. Pasechnik, On the irreducibility of a truncated
binomial expansion,\emph{ Rocky Mountain J. Math.} 37 (2007), 455-464.\\
{[3] I. V. Ostrovskii, On a problem of A. Eremenko, \emph{Comput. Methods Funct. Theory} 4 (2004), No. 2,
275-282.\\
{[4] P. Ribenboim, \textit{The Theory of Classical Valuations}, Springer-$\textrm{Verlag}$ New York, 1999. \\
{[5] S. Laishram and T. N. Shorey, Number of prime divisors in a product of
consecutive integers, \emph{Acta Arithmetica} 113.4 (2004), 327-341.\\
{[6] I. Scherbak, Intersections of Schubert varieties and highest weight vectors in tensor products $sl_{N+1^-}$
representations, ArXiv e-print math.RT/0409329, Sept. 2004.\\
{[7] N. Saradha and T. N. Shorey, Almost squares and factorisations in consecutive
integers, \emph{Compositio Math}. 138 (2003), 113-124.\\
\vspace{3mm}\\

\baselineskip 10pt
\noindent Sudesh K. Khanduja, Ramneek Khassa~~~~~~~Shanta Laishram\\
Department of Mathematics,~~~~~~~~~~~~~~~~~~~~Stat-Math Unit, Indian Statistical Institute,\\
Panjab University, Chandigarh-160014, ~~~~~~~~7, S.J.S. Sansanwal Marg, New Delhi-110016. \\
India.~~~~~~~~~~~~~~~~~~~~~~~~~~~~~~~~~~~~~~~~~~~~~~~~~~~India.\\
Email: skhand@pu.ac.in, ramneekkhassa@yahoo.co.in, shanta@isid.ac.in

\end{document}